\documentclass[10pt]{amsart}
\usepackage{amsfonts}
\usepackage{color}
\usepackage[square,compress,comma, numbers]{natbib}
\usepackage[colorlinks=true, citecolor=blue, linkcolor=blue]{hyperref}
\allowdisplaybreaks[4]
\usepackage{amssymb}
\usepackage{amsmath}
\usepackage{color}

\definecolor{c20}{rgb}{0.,0.7,0.}
\definecolor{c30}{rgb}{0.,0.,1.}
\definecolor{c40}{rgb}{1,0.1,0.7}
\definecolor{c50}{rgb}{1,0,0}
\definecolor{c60}{rgb}{1,0.9,0.1}

\allowdisplaybreaks[4]

\newcommand{\abs}[1]{\left\lvert #1 \right\rvert}

\newcommand{\E}[1]{\mathbb{E}\left\{#1\right\}}

\newcommand{\pk}[1]{\mathbb{P} \left \{#1 \right \} }
\newcommand{\R}{\mathbb{R}}
\newcommand{\N}{\mathbb{N}}

\newcommand{\BQN}{\begin{eqnarray}}
\newcommand{\EQN}{\end{eqnarray}}
\newcommand{\BQNY}{\begin{eqnarray*}}
\newcommand{\EQNY}{\end{eqnarray*}}

\newcommand{\BS}{\begin{sat}}
\newcommand{\ES}{\end{sat}}
\newcommand{\BT}{\begin{theo}}
\newcommand{\ET}{\end{theo}}
\newcommand{\BL}{\begin{lem}}
\newcommand{\EL}{\end{lem}}
\newcommand{\BK}{\begin{korr}}
\newcommand{\EK}{\end{korr}}

\newcommand{\BD}{\begin{de}}
\newcommand{\ED}{\end{de}}
\newcommand{\BIT}{\begin{itemize}}
\newcommand{\EIT}{\end{itemize}}
\newcommand{\BDI}{\begin{description}}
\newcommand{\EDI}{\end{description}}

\newcommand{\BRM}{\begin{remarks}}
\newcommand{\ERM}{\end{remarks}}

\newcommand{\BEL}{\begin{lem}}
\newcommand{\EEL}{\end{lem}}

\newtheorem{theo}{Theorem}[section]
\newtheorem{sat}[theo]{Proposition}
\newtheorem{de}[theo]{Definition}
\newtheorem{lem}[theo]{Lemma}

\newtheorem{example}[theo]{Example}
\newtheorem{korr}[theo]{Corollary}
\newtheorem{remark}[theo]{Remark}
\newtheorem{remarks}[theo]{Remarks}

\newcommand{\nelem}[1]{{Lemma \ref{#1}}}

\newcommand{\netheo}[1]{{Theorem \ref{#1}}}

\newcommand{\prooftheo}[1]{ \textsc{\bf Proof of Theorem} \ref{#1}:}

\newcommand{\prooflem}[1]{\textsc{\bf Proof of Lemma} \ref{#1}:}

\newcommand{\COM}[1]{}

\newcommand{\QED}{\hfill $\Box$}

\def\en{\epsilon}

\def\rw{\rightarrow}

\def\IF{\infty}

\def\LT{\left}
\def\RT{\right}

\def\rw{\rightarrow}

\def\vn{\varepsilon}
\def\Var{\text{Var}}

\def\Bu+#1{\mathcal{B}^{\varepsilon+}_{u}(#1)}

\begin{document}

\title{Extremes of $\alpha(t)$-Locally stationary Gaussian processes with non-constant variances}

\author{Long Bai}
\address{Long Bai, Department of Actuarial Science, 
University of Lausanne\\
UNIL-Dorigny, 1015 Lausanne, Switzerland
}
\email{Long.Bai@unil.ch}
\bigskip

 \maketitle

{\bf Abstract:} With motivation from \cite{atlocally}, in this paper we derive the exact tail asymptotics of $\alpha(t)$-locally stationary
Gaussian processes with non-constant variance functions. We show that some certain variance functions lead to qualitatively new results.


\section{Introduction and Main Result }
For $X(t),\ t\in [0,T],\ T>0$ a centered stationary Gaussian process with unit variance and continuous sample paths Pickands derived in \cite{PicandsA} that
\begin{eqnarray}\label{eq1.2}
\pk{\sup_{t\in [0,T]} X(t)>u}\sim T\mathcal{H}_{\alpha}a^{1/\alpha}u^{2/\alpha}\pk{X(0)> u}, \quad  \ u\rightarrow\infty,
\end{eqnarray}
provided that the correlation function $r$ satisfies
\BQN \label{PickandsC}
1-r(t) \sim  a\abs{t}^{\alpha} , \quad t \downarrow 0,  \quad a>0,  \quad \text{ and  } r(t)< 1, \ \forall\ t\not=0,
\EQN
with $\alpha \in (0,2]$ ($\sim$ means asymptotic equivalence when the argument tends to 0 or $\IF$). Here the classical Pickands constant $\mathcal{H}_{\alpha}$ is defined by
$$\mathcal{H}_{\alpha}=\lim_{T\rightarrow\infty} T^{-1} \E{ \sup_{t\in[0,T]} e^{\sqrt{2}B_{\alpha}(t)-t^{\alpha}}},$$
where ${B_{\alpha}(t),t\ge 0}$  is a standard fractional Brownian motion with Hurst index $\alpha/2\in (0,1]$, see
 \cite{PicandsA,Pit72, DE2002,DI2005,DE2014,DiekerY,DEJ14,Pit20, Tabis, DM, SBK, Htilt} for various properties of $\mathcal{H}_\alpha$.\\
The deep contribution \cite{Berman92} introduced the class of locally stationary Gaussian processes with index $\alpha$, i.e., a centered Gaussian process $X(t),t\in[0,T]$ with a constant variance function, say equal to 1, and correlation function satisfying $$r(t,t+h)=1-a(t)|h|^\alpha+o(|t|^\alpha), \ h\rw 0,$$ uniformly with respect to $t\in[0,T]$, where $\alpha\in(0,2]$ and $a(t)$ is a bounded, strictly positive and continuous function.\\
Clearly, the class of locally stationary Gaussian processes includes the stationary ones. It allows for some minor fluctuations of dependence at $t$ and at the same time keeps stationary structure at the local scale. See \cite{Berman92,braker1993high,Husler90} for studies on the  locally stationary Gaussian processes with index $\alpha$.

In \cite{atlocally} the tail asymptotics of the supremum of $\alpha(t)$-locally stationary Gaussian processes  are investigated. Such processes and random fields  are of interest in various applications, see \cite{atlocally} and the recent contributions \cite{Xiao,MR3413855,Lif}. Following the definition in \cite{atlocally}, a centered Gaussian process $X(t),t\in[0,T]$ with  continuous sample paths and unit variance is $\alpha(t)$-locally stationary if the correlation function $r(\cdot,\cdot)$ satisfies the following conditions:
\begin{itemize}
\item [(i)] $\alpha(t)\in C([0,T])$ and $\alpha(t)\in(0,2]$ for all $t\in[0,T]$;
\item [(ii)] $a(t)\in C([0,T])$ and $0<\inf\{a(t):t\in[0,T]\}\leq\sup\{a(t):t\in[0,T]\}<\IF$;
\item [(iii)]  uniformly for $t\in[0,T]$
$$1-r(t,t+h)=a(t)|h|^{\alpha(t)}+o(|h|^{\alpha(t)}),\ h\rw 0,$$
where $f(t)\in C(\mathcal{T})$ means that $f(t)$ is continuous on $\mathcal{T}\subset\R$.
\end{itemize}

In this paper, we shall consider the case that the variance function $\sigma^2(t)=Var(X(t))$ is not constant, assuming instead that:  \\
(iv) $\sigma(t)$ attains its maximum equal to 1 over $[0,T]$ at the unique point $t_0\in[0,T]$ and for some constants $c, \gamma>0$,  $$\frac{1}{\sigma(t)}=1+ce^{-|t-t_0|^{-\gamma}}(1+o(1)), \ \ t\rw t_0.$$
A crucial assumption in our result is that similar  to the variance function, the function  $\alpha(t)$ has a certain behaviour around the extreme point $t_0$.
Specifically,  as in \cite{atlocally} we shall assume: \\
(v) there exist $ \beta, \delta,\ b>0$ such that
$$\alpha(t+t_0)=\alpha(t_0)+b|t|^\beta+o(|t|^{\beta+\delta}),\quad t\rw 0.$$
\begin{remark}
We remark that $t_0$ does not need to be the unique point such that $\alpha(t)$ is minimal on $[0,T]$, which is different from  \cite{atlocally}. For instance, $[0,T]=[0,2\pi]$, $t_0=0$ and $\alpha(t)=1+\frac{1}{2}\sin(t)$, then $0$ is not the minimum point of $\alpha(t)$ over $[0,2\pi]$ which means assumptions about $\alpha(t)$ in \cite{atlocally} are not satisfied but assumption (v) here is satisfied with
$$
\alpha(t)=1+\frac{1}{2}|t|+o(|t|^{\frac{3}{2}}),\ t\rw 0.
$$
\end{remark}
Below we set  $\alpha:=\alpha(t_0)$, $a:=a(t_0)$ and write $\Psi$ for the survival function of an $N(0,1)$ random variable.
Further, define $0^{a} =\IF$ for $a<0$. Our main result is stated in the next theorem.
\BT\label{Main1}
If a centered Gaussian process $X(t), t\in [0,T]$ with continuous sample paths is such that the assumptions \hbox{(i)-(v)} are valid, then we have as $u\to \IF$
\BQNY
\pk{\sup_{t\in[0,T]}X(t)>u}\sim\widehat{I}a^{1/\alpha}\mathcal{H}_{\alpha}u^{2/\alpha}(\ln u)^{-\frac{1}{\gamma\wedge\beta}}\Psi(u)
\LT\{
\begin{array}{ll}
	2 ^{-1/\gamma},& \hbox{if}\ \gamma<\beta,\\
	\int_{0}^{2^{-1/\gamma}}e^{\frac{-2b x^\beta}{\alpha^2}}dx, & \hbox{if}\ \gamma=\beta,\\
	\int_{0}^{\IF}e^{\frac{-2b x^\beta}{\alpha^2}}dx, & \hbox{if}\ \gamma>\beta,
\end{array}
\RT.
\EQNY
where $\gamma\wedge\beta=\min (\gamma, \beta)$ and
\BQNY
\widehat{I}=
\LT\{
\begin{array}{ll}
	1, & \hbox{if}\ t_0=0\ \hbox{or}\ t_0=T,\\
	2, & \hbox{if}\ t_0\in(0,T).
\end{array}
\RT.
\EQNY
\ET
\begin{remark}
	i) If $\alpha(t)\equiv \alpha$ for all t in a small neighborhood of $t_0$, the asymptotic of $\pk{\sup_{t\in[0,T]}X(t)>u}$ is the same as in the case of $\gamma<\beta$ in \netheo{Main1}.\\
	ii) The result of case $\gamma>\beta$ in \netheo{Main1} is the same as the $\alpha(t)$-locally stationary scenario in \cite{atlocally}, which means that $\sigma(t)$ varies so slow in a small neighborhood of $t_0$ that $X(t)$ can be considered as $\alpha(t)$-locally stationary in this small neighborhood.
\end{remark}
The following example is a straightforward application of \netheo{Main1}.
\begin{example}
Here we consider a multifractional Brownian motion  $B_{H(t)}(t), \ t\geq 0$, i.e., a centered Gaussian process with covariance function
\BQNY
\E{B_{H(t)}(t)B_{H(s)}(s)}=\frac{1}{2}D(H(s)+H(t))\LT[|s|^{H(s)+H(t)}+|t|^{H(s)+H(t)}-|t-s|^{H(s)+H(t)}\RT],
\EQNY
where $D(x)=\frac{2\pi}{\Gamma(x+1)\sin\LT(\frac{\pi x}{2}\RT)}$ and $H(t)$ is a H\"{o}lder function of exponent $\lambda$ such that $0<H(t)<\min(1,\lambda)$ for $t\in[0,\IF)$. For constants $T_1,T_2$ with $0<T_1<T_2$, define
\BQNY
\overline{B}_{H(t)}(t):=\frac{B_{H(t)}(t)}{\sqrt{\Var(B_{H(t)}(t))}},\ t\in[T_1,T_2],
\EQNY
and
\BQNY
\sigma(t):=1-e^{-|t-t_0|^{-\gamma}}, \ t\in[T_1,T_2],
\EQNY
with some $t_0\in(T_1,T_2)$ and $\gamma>0$.\\
By \cite{atlocally}, $\overline{B}_{H(t)}(t),\ t\in[T_1,T_2],$ is a $2H(t)$-locally stationary Gussian process with correlation function 
\BQNY
r(t,t+h)=1-\frac{1}{2}t^{-2H(t)}|h|^{2H(t)}+o(|h|^{2H(t)}),\ \ h\rw 0.
\EQNY
Further, we assume that there exist $\beta, \delta, b>0$ such that $H(t+t_0)=H(t_0)+bt^{\beta}+o(t^{\beta+\delta})$, as $t\rw 0$. Then
\BQNY
\pk{\sup_{t\in[T_1,T_2]}\sigma(t)\overline{B}_{H(t)}(t)>u}\sim
2^{1-1/2H}\frac{\mathcal{H}_{2H}}{t_0}u^{1/H}(\ln u)^{-\frac{1}{\gamma\wedge \beta}}\Psi(u)
\LT\{
\begin{array}{ll}
	2 ^{-1/\gamma},& \hbox{if}\ \gamma<\beta,\\
	\int_{0}^{2^{-1/\gamma}}e^{\frac{-b x^\beta}{H^2}}dx, & \hbox{if}\ \gamma=\beta,\\
	\int_{0}^{\IF}e^{\frac{-b x^\beta}{H^2}}dx, & \hbox{if}\ \gamma>\beta,
\end{array}
\RT.
\ \ u\rw 0.
\EQNY
with $H:=H(t_0)$.
\end{example}

\section{Proofs}
In the rest of the paper, we focus on the case when $t_0=0$. The complementary scenario when $t_0\in(0,T]$ follows by analogous argumentation. 
 Recall that
$$\mathcal{H}_{\alpha}=\lim_{T\rightarrow\IF}\frac{1}{T}\mathcal{H}_{\alpha}[0,T],\quad \text{with }
\mathcal{H}_{\alpha}[-S_1,S_2]=\E{\sup_{t\in[-S_1,S_2]}e^{\sqrt{2}B_\alpha(t)-|t|^\alpha}}\in(0,\IF),
$$
where $S_1,S_2 \in [0,\IF)$ with $\max(S_1,S_2)>0$ are some constants.

\BEL\label{le1}
Under the assumptions of \netheo{Main1} we have
\BQN\label{le1eq1}
\pk{\sup_{t\in[0,T]}X(t)>u}\sim\pk{\sup_{t\in[0,\delta_1(u)]}X(t)>u}, \ u\rw\IF.
\EQN
Moreover, there exists a constant $C>0$ such that for all sufficiently large $u$
\BQN\label{le1eq2}
\pk{\sup_{t\in[\delta_2(u),T]}X(t)>u}\leq C T u^{2/\alpha} (\ln u)^{-4/3\beta}\Psi\LT(u\RT),
\EQN
where for some constant $q>1$
\BQN\label{de}
\delta_1(u)=\LT(\frac{1}{2\ln u-q\ln\ln u}\RT)^{1/\gamma}\ \ \text{and}\ \ \delta_2(u)=\LT(\frac{\alpha^2 (\ln (\ln u))}{\beta(\ln u)}\RT)^{1/\beta}.
\EQN
\EEL
\vskip 0.5cm
By \eqref{le1eq2}, in the proof of \netheo{Main1}, we derive that, as $u\rw\IF$,
\BQN\label{ne}
\pk{\sup_{t\in[\delta_2(u),T]}X(t)>u}=o\LT(\pk{\sup_{t\in[0,\delta_2(u)]}X(t)>u}\RT).
\EQN
Since $\delta_1(u)\rw 0, \delta_2(u)\rw 0$ as $u\rw \IF$ and $a(t)$ is continuous,  without loss of generality, we may assume that $a(t)\equiv a(0)=a$ for $t\in\LT([0,\delta_1(u)]\cup[0,\delta_2(u)]\RT)$.
Moreover, by assumption (iv), we know that $\sigma(t)> 0$ for $t \in \LT([0,\delta_1(u)]\cup[0,\delta_2(u)]\RT)$. Below we use notation $\overline{X}(t)=\frac{X(t)}{\sigma(t)}$ for all $t$ such  that $\sigma(t)$ is positive.
\vskip 0.5cm
\prooftheo{Main1}
 First we derive the asymptotic of
\BQNY
\pi(u):=\pk{\sup_{t\in\Delta(u)}X(t)>u},
\EQNY
as $u\rw\IF$, where $\Delta(u)=[0,\delta(u)]$ and
\BQNY
\delta(u)=
\LT\{
\begin{array}{ll}
\delta_1(u), & \hbox{if}\ \gamma\leq \beta,\\
\delta_2(u), & \hbox{if}\ \gamma>\beta,
\end{array}
\RT.
\EQNY
with $\delta_1(u)$ and $\delta_2(u)$  in \eqref{de}, which combined with \nelem{le1} finally shows that
\BQN\label{nee}
\pk{\sup_{t\in[0,T]}X(t)>u}\sim\pi(u).
\EQN
In the following  $\mathbb{Q}_i, \ i\in\N,$ are some positive constants.
For some $S>0$, let
$Y_{\nu,u}(t), t\in[0,S]$ be a family of centered stationary Gaussian processes with
\BQNY
Cov\LT(Y_{\nu,u}(s),Y_{\nu,u}(t)\RT)=1-(1-\nu)au^{-2}|s-t|^{\alpha+2b\delta^\beta(u)},
\EQNY
for $\nu\in(0,1), u>0$ such that $\alpha+2b\delta^\beta(u)\leq 2$ and $s,t\in[0,S]$.
Further, let  $Z_{\nu,u}(t), t\in[0,S]$ be another  family of centered stationary Gaussian processes with
\BQNY
Cov\LT(Z_{\nu,u}(s),Z_{\nu,u}(t)\RT)=1-(1+\nu)au^{-2}|s-t|^{\alpha},
\EQNY
for $\nu\in(0,1), u>0$ and $s,t\in[0,S]$.
Due to assumptions (i) and (v), $\alpha$ is strictly smaller than 2, which guarantees that covariance function of $Y_{\nu,u}(t), t\in[0,S]$ and $Z_{\nu,u}(t), t\in[0,S]$ are positive-definite. Hence the introduced families of Gaussian processes exist.\\
By assumption (iv), for any small $\vn\in(0,1)$
\BQN\label{sig}
1+(1-\vn)ce^{-|t|^{-\gamma}}\leq\frac{1}{\sigma(t)}\leq1+(1+\vn)ce^{-|t|^{-\gamma}},
\EQN
holds for $t\in [0,\delta(u)]$.\\
{\bf Case 1}: $\gamma<\beta$.
Set for any $\en \in(0,1)$ and all $u$ large
\BQNY
&& N(0)=N(u,0):=\LT\lfloor\frac{\delta_1(u)u^{2/\alpha}}{S}\RT\rfloor,\ \
 N_{\en}(u)=\LT\lfloor(1-\en)\frac{\delta_1(u)u^{2/\alpha}}{S}\RT\rfloor=\LT\lfloor\frac{(1-\en)u^{2/\alpha}}{(2\ln u-q\ln\ln u)^{1/\gamma}S}\RT\rfloor,\\
&&B_j(u)=B_{j,0}(u)=\LT[j\frac{S}{u^{2/\alpha}},(j+1)\frac{S}{u^{2/\alpha}}\RT],\ j\in \N, \ \
\mathcal{G}_u^{\pm\vn}=u\LT(1+(1\pm\vn)ce^{-((1-\en)\delta_1(u))^{-\gamma}}\RT).
\EQNY
We notice the fact that
\BQNY
\Psi(\mathcal{G}_u^{\pm\vn})\sim\Psi(u),\ u\rw\IF,
\EQNY
and
\BQN\label{bound1}
I_1(u)\leq\pi(u)\leq I_1(u)+I_2(u),
\EQN
where
\BQNY
I_1(u)&=&\pk{\sup_{t\in[0,(1-\en)\delta_1(u)]}X(t)>u}, \ \ I_2(u)=\pk{\sup_{t\in[(1-\en)\delta_1(u),\delta_1(u)]}X(t)>u}.
\EQNY
Then by Bonferroni's inequality, \eqref{sig}, \nelem{lein} with $k=0$ and \nelem{leYZ}
\BQN\label{I1upper1}
I_1(u)
&\leq&\sum_{j=0}^{N_{\en}(u)}\pk{\sup_{t\in B_j(u)}X(t)>u}\nonumber\\
&\leq&\sum_{j=0}^{N_{\en}(u)}\pk{\sup_{t\in B_j(u)}\overline{X}(t)>\mathcal{G}_u^{-\vn}}\nonumber\\
&\leq&\sum_{j=0}^{N_{\en}(u)}\pk{\sup_{t\in [jS,(j+1)S]}\overline{X}(tu^{-2/\alpha})>\mathcal{G}_u^{-\vn}}\nonumber\\
&\leq&\sum_{j=0}^{N_{\en}(u)}\pk{\sup_{t\in [0,S]}Z_{u,\nu}(t)>\mathcal{G}_u^{-\vn}}\nonumber\\
&\sim&\sum_{j=0}^{N_{\en}(u)}\mathcal{H}_{\alpha}\LT[0,S((1+\nu)a)^{1/\alpha}\RT]\Psi\LT(\mathcal{G}_u^{-\vn}\RT)\nonumber\\
&\sim&\sum_{j=0}^{N_{\en}(u)}\mathcal{H}_{\alpha}\LT[0,S((1+\nu)a)^{1/\alpha}\RT]\Psi(u)\nonumber\\
&\sim&(1-\en)u^{2/\alpha}\delta_1(u)\frac{\mathcal{H}_{\alpha}\LT[0,S((1+\nu)a)^{1/\alpha}\RT]}{S}\Psi(u)\nonumber\\
&\sim&(1-\en)((1+\nu)a)^{1/\alpha}\mathcal{H}_{\alpha}u^{2/\alpha}\delta_1(u)\Psi(u), \ u\rw\IF, \ S\rw\IF.
\EQN
Similarly,
\BQN\label{lower1}
\sum_{j=0}^{N_{\en}(u)-1}\pk{\sup_{t\in B_j(u)}X(t)>u}&\geq&\sum_{j=0}^{N_{\en}(u)-1}\pk{\sup_{t\in [0,S]}Y_{u,\nu}(t)>\mathcal{G}_u^{+\vn}}\nonumber\\
&\sim&(1-\en)((1-\nu)a)^{1/\alpha}\mathcal{H}_{\alpha}u^{2/\alpha}\delta_1(u)\Psi(u), \ u\rw\IF, \ S\rw\IF.
\EQN
Since
\BQN\label{I1lower1}
I_1(u)\geq \sum_{j=0}^{N_{\en}(u)-1}\pk{\sup_{t\in B_j(u)}X(t)>u}
-\sum_{0\leq j<k\leq N_{\en}(u)}\pk{\sup_{t\in B_j(u)}X(t)>u,\sup_{t\in B_k(u)}X(t)>u},
\EQN
and by \cite{atlocally}[Lemma 4.5]
\BQN\label{lower2}
\sum_{0\leq j<k\leq N_{\en}(u)}\pk{\sup_{t\in B_j(u)}X(t)>u,\sup_{t\in B_k(u)}X(t)>u}
&\leq&\sum_{0\leq j<k\leq N_{\en}(u)}\pk{\sup_{t\in B_j(u)}\overline{X}(t)>u,\sup_{t\in B_k(u)}\overline{X}(t)>u}\nonumber\\
&=&o\LT(u^{2/\alpha}\delta_1(u)\Psi(u)\RT),\ u\rw\IF,\ S\rw\IF, \ \en\rw 0.
\EQN
Thus inserting \eqref{lower1} and \eqref{lower2} into \eqref{I1lower1}, we have
\BQNY
\lim_{u\rw\IF}\frac{I_1(u)(2\ln u-q\ln\ln u)^{1/\gamma}}{u^{2/\alpha}\Psi(u)}\geq (1-\en)((1-\nu)a)^{1/\alpha}\mathcal{H}_{\alpha},
\EQNY
which combined with \eqref{I1upper1} gives that
\BQN\label{1I1}
I_1(u)\sim \frac{a^{1/\alpha}\mathcal{H}_{\alpha}u^{2/\alpha}}{(2\ln u-q\ln\ln u)^{1/\gamma}}\Psi(u), u\rw\IF,\ \nu\rw 0,\ \en\rw0.
\EQN

By (iii) and (v), we have for all $u$ large
\BQNY
\E{(\overline{X}(t)-\overline{X}(s))^2}=2-2r(s,t)\leq
\mathbb{Q}_1|s-t|^{\alpha},
\EQNY
uniformly holds for $s,t\in[(1-\en)\delta_1(u),\delta_1(u)]$.
By Piterbarg inequality for $u$ large enough, see e.g., \cite{Pit96}[Theorem 8.1] or an extension in \cite{KEP2015}[Lemma 5.1]
\BQN\label{I2upper}
I_2(u)\leq\pk{\sup_{t\in[(1-\en)\delta_1(u),\delta_1(u)]}\overline{X}(t)>u}\leq
\mathbb{Q}_2\en\delta_1(u)u^{2/\alpha}\Psi(u),
\EQN
which implies
\BQNY
\lim_{\en\rw 0}\lim_{u\rw\IF}\frac{I_2(u)(2\ln u-q\ln\ln u)^{1/\gamma}}{u^{2/\alpha}\Psi(u)}=0.
\EQNY
Combining this equation with \eqref{bound1} and \eqref{1I1}, we get
$$\pi(u)\sim \frac{a^{1/\alpha}\mathcal{H}_{\alpha}u^{2/\alpha}}{(2\ln u-q\ln\ln u)^{1/\gamma}}\Psi(u)
\sim a^{1/\alpha}\mathcal{H}_{\alpha}u^{2/\alpha}(2\ln u)^{-1/\gamma}\Psi(u), \ u\rw\IF.$$
{\bf Case 2}: $\gamma=\beta$.
Set
\BQNY
d_k=d_k(u):=\LT(\frac{k}{\ln (u)(\ln \ln (u))^{1/\beta}}\RT)^{1/\beta},\ \ A_k=A_k(u):=\LT[d_k,d_{k+1}\RT].
\EQNY
Further let $M_{\en}(u)=\max(k\in\N: d_k\leq (1-\en)\delta_1(u))$ for some $\en\in(0,1)$, then $M_{\en}(u)\rw \IF, \ u\rw\IF$. Clearly
$$\bigcup_{k=0}^{M_{\en}(u)-1}A_k\subset[0,(1-\en)\delta_1(u)]\subset\bigcup_{k=0}^{M_{\en}(u)}A_k.$$
We divide each interval $A_k$ into subintervals of length $S/u^{2/\alpha(d_{k})}$, i.e.,
\BQNY
B_{j,k}=B_{j,k}(u):=\LT[d_k+j\frac{S}{u^{2/\alpha(d_{k})}},d_k+(j+1)\frac{S}{u^{2/\alpha(d_{k})}}\RT]
\EQNY
for $j=0,1,\ldots, N(k)$, where $N(k)=N(k,u):=\LT\lfloor\frac{d_{k+1}-d_k}{S}u^{2/\alpha(d_{k})}\RT\rfloor$. Notice that
$$\bigcup_{k=0}^{N(k)-1}B_{j,k}\subset A_k\subset\bigcup_{k=0}^{N(k)}B_{j,k}.$$
We have
\BQN\label{bound2}
I_1(u)\leq\pi(u)\leq I_1(u)+I_2(u),
\EQN
where
\BQNY
I_1(u)&=&\pk{\sup_{t\in[0,(1-\en)\delta_1(u)]}X(t)>u}, \ \ I_2(u)=\pk{\sup_{t\in[(1-\en)\delta_1(u),\delta_1(u)]}X(t)>u}.
\EQNY
Then by Bonferroni's inequality
\BQN\label{boundlower}
I_1(u)&\geq&\sum_{k=0}^{M_{\en}(u)-1}\sum_{j=0}^{N(k)-1}\pk{\sup_{t\in B_{j,k}}X(t)>u}-
\underset{(j,k)\prec(j',k')}{\sum_{(j,k),(j',k')\in\mathcal{L}}}\pk{\sup_{t\in B_{j,k}}X(t)>u,\sup_{t\in B_{j',k'}}X(t)>u}\nonumber\\
&=:&J_1(u)-J_2(u),
\EQN
where $\mathcal{L}=\{(j,k):0\leq k\leq M_{\en}(u)-1, 0\leq j\leq N(k)-1\}$ and
$$(j,k)\prec(j',k')\ \hbox{iff}\ (k<k')\vee(k=k'\wedge j<j'),$$
and by \eqref{sig}, \nelem{lein} and \nelem{leYZ}
\BQNY
I_1(u)&\leq&\sum_{k=0}^{M_{\en}(u)}\sum_{j=0}^{N(k)}\pk{\sup_{t\in B_{j,k}}X(t)>u}\\
&\leq&\sum_{k=0}^{M_{\en}(u)}\sum_{j=0}^{N(k)}\pk{\sup_{t\in B_{j,k}}\overline{X}(t)>\mathcal{G}_u^{-\vn}}\\
&\leq&\sum_{k=0}^{M_{\en}(u)}\sum_{j=0}^{N(k)}\pk{\sup_{t\in [0,S]}Z_{\nu,u}(t)>\mathcal{G}_u^{-\vn}}\\
&\sim&\sum_{k=0}^{M_{\en}(u)}\sum_{j=0}^{N(k)}\mathcal{H}_{\alpha}\LT[0,S((1+\nu)a)^{1/\alpha}\RT]\Psi\LT(\mathcal{G}_u^{-\vn}\RT)\\
&\sim&\sum_{k=0}^{M_{\en}(u)}\frac{d_{k+1}-d_k}{S}u^{2/\alpha(d_{k})}\mathcal{H}_{\alpha}\LT[0,S((1+\nu)a)^{1/\alpha}\RT]\Psi(u)\\
&=&\frac{\mathcal{H}_{\alpha}\LT[0,S((1+\nu)a)^{1/\alpha}\RT]}{S}\frac{u^{2/\alpha}}{(\ln u)^{1/\beta}}\Psi(u)\sum_{k=0}^{M_{\en}(u)}(\ln u)^{1/\beta}(d_{k+1}-d_k)e^{(\ln u)\LT(\frac{2(\alpha-\alpha(d_{k}))}{\alpha\alpha(d_{k})}\RT)}\\
&\leq&\frac{\mathcal{H}_{\alpha}\LT[0,S((1+\nu)a)^{1/\alpha}\RT]}{S}\frac{u^{2/\alpha}}{(\ln u)^{1/\beta}}\Psi(u)\sum_{k=0}^{M_{\en}(u)}(\ln u)^{1/\beta}(d_{k+1}-d_k)e^{\frac{-2(1-\vn_1)(\ln u)\LT(bd^\beta_{k}-d^{\beta+\delta}_{k}\RT)}{\alpha^2}}\\
&\leq&\frac{\mathcal{H}_{\alpha}\LT[0,S((1+\nu)a)^{1/\alpha}\RT]}{S}\frac{u^{2/\alpha}}{(\ln u)^{1/\beta}}\Psi(u)\\
&&\times\sum_{k=0}^{M_{\en}(u)}(\ln u)^{1/\beta}(d_{k+1}-d_k)e^{\frac{-2(1-\vn_1)b\LT((\ln u)^{1/\beta}d_{k}\RT)^\beta}{\alpha^2}}e^{\frac{2(1-\vn_1)(\ln u)d^{\beta+\delta}_{M_{\en}(u)+1}}{\alpha^2}},
\EQNY
as $u\rw\IF$, where $\vn_1\in (0,1)$ is a small constant.\\
Moreover, using that $d_{M_{\en}(u)}\leq(1-\en)\delta_1(u)$ and
$\lim_{u\rw\IF}(\ln u)\delta_1(u)^{\beta+\delta}=0$, we observe that
\BQNY
\lim_{u\rw\IF}e^{\frac{2(1-\vn_1)(\ln u)d^{\beta+\delta}_{M_{\en}(u)+1}}{\alpha^2}}=1.
\EQNY
Finally, since $$\lim_{u\rw\IF}\sup_{k=0,\ldots,M_{\en}(u)}(\ln u)^{1/\beta}(d_{k+1}-d_{k})=0$$  and
$$\lim_{u\rw\IF}(\ln u)^{1/\beta}d_{M_{\en}(u)+1}= (1-\en)\LT(\frac{1}{2}\RT)^{1/\beta},$$
we obtain
\BQNY
\lim_{u\rw\IF}\sum_{k=0}^{M_{\en}(u)}(\ln u)^{1/\beta}(d_{k+1}-d_k)e^{\frac{-2(1-\vn_1)b\LT((\ln u)^{1/\beta}d_{k}\RT)^\beta}{\alpha^2}}=\int_{0}^{(1-\en)\LT(\frac{1}{2}\RT)^{1/\beta}}e^{\frac{-2(1-\vn_1)b x^\beta}{\alpha^2}}dx.
\EQNY
Thus
\BQN\label{2I1upper}
\lim_{u\rw\IF}\frac{I_1(u)(\ln u)^{1/\beta}}{u^{2/\alpha}\Psi(u)}\leq\frac{\mathcal{H}_{\alpha}\LT[0,S((1+\nu)a)^{1/\alpha}\RT]}{S}\int_{0}^{(1-\en)\LT(\frac{1}{2}\RT)^{1/\beta}}e^{\frac{-2(1-\vn_1)b x^\beta}{\alpha^2}}dx,
\EQN
and letting $S\rw \IF,\vn_1,\nu \rw 0$, and $\en\rw0$, we get the upper bound.
Similarly, we derive that
\BQN\label{J1}
\lim_{\en\rw 0}\lim_{S\rw\IF}\lim_{u\rw\IF}\frac{J_1(u)(\ln u)^{1/\beta}}{u^{2/\alpha}\Psi(u)}
\geq a^{1/\alpha}\mathcal{H}_{\alpha}\int_{0}^{\LT(\frac{1}{2}\RT)^{1/\beta}}e^{\frac{-2b x^\beta}{\alpha^2}}dx.
\EQN
By \cite{atlocally} [Lemma 4.5]
\BQN\label{lower3}
J_2(u)&=&\underset{(j,k)\prec(j',k')}{\sum_{(j,k),(j',k')\in\mathcal{L}}}\pk{\sup_{t\in B_{j,k}}X(t)>u,\sup_{t\in B_{j',k'}}X(t)>u}\nonumber\\
&\leq&\underset{(j,k)\prec(j',k')}{\sum_{(j,k),(j',k')\in\mathcal{L}}}\pk{\sup_{t\in B_{j,k}}\overline{X}(t)>u,\sup_{t\in B_{j',k'}}\overline{X}(t)>u}\nonumber\\
&=&o\LT(u^{2/\alpha}(\ln u)^{-1/\beta}\Psi(u)\RT),\ u\rw\IF,\ S\rw\IF, \en\rw 0.
\EQN
Thus inserting \eqref{J1} and \eqref{lower3} into \eqref{boundlower}, we get
\BQN\label{2I1lower}
\lim_{\en\rw 0}\lim_{S\rw\IF}\lim_{u\rw\IF}\frac{I_1(u)(\ln u)^{1/\beta}}{u^{2/\alpha}\Psi(u)}\geq a^{1/\alpha}\mathcal{H}_{\alpha}\int_{0}^{\LT(\frac{1}{2}\RT)^{1/\beta}}e^{\frac{-2(1-\vn_1)b x^\beta}{\alpha^2}}dx.
\EQN
By \eqref{I2upper}
\BQN\label{2I2upper}
\lim_{\en\rw 0}\lim_{u\rw\IF}\frac{I_2(u)(\ln u)^{1/\beta}}{u^{2/\alpha}\Psi(u)}=0.
\EQN
Hence according to \eqref{bound2}, \eqref{2I1upper}, \eqref{2I1lower},  and  \eqref{2I2upper}, we have
$$\pi(u)\sim a^{1/\alpha}\mathcal{H}_{\alpha}u^{2/\alpha}(\ln u)^{-1/\beta}\Psi(u)\int_{0}^{\LT(\frac{1}{2}\RT)^{1/\beta}}e^{\frac{-2b x^\beta}{\alpha^2}}dx, \ u\rw\IF.$$
{\bf Case 3}: $\gamma>\beta$.
We consider $\pi(u)=\pk{\sup_{t\in[0,\delta_2(u)]}X(t)>u}$ with
\BQNY
\delta_2(u)=\LT(\frac{\alpha^2 (\ln (\ln u))}{\beta(\ln u)}\RT)^{1/\beta}.
\EQNY
Set for some $\vn>0$
\BQNY
\mathcal{F}_u^{\pm\vn}=u\LT(1+(1\pm\vn)ce^{-(\delta_2(u))^{-\gamma}}\RT),\
\mathcal{K}=\{t\in[0,T]:\sigma(t)\neq 0\},
\EQNY
and we observe that  $$\Psi\LT(\mathcal{F}_u^{\pm\vn}\RT)\sim\Psi(u), \ u\rw\IF.$$
By \cite{atlocally}[Theorem 2.1]
\BQN\label{3upper}
\pi(u)&\leq&\pk{\sup_{t\in[0,\delta_2(u)]}\overline{X}(t)>u}\nonumber\\
&\leq&\pk{\sup_{t\in\mathcal{K}}\overline{X}(t)>u}\nonumber\\
&\sim&a^{1/\alpha}\mathcal{H}_{\alpha}u^{2/\alpha}(\ln u)^{-\frac{1}{\beta}}\int_{0}^{\IF}e^{\frac{-2b x^\beta}{\alpha^2}}dx\Psi(u), \ u\rw\IF.
\EQN
Let $d_k, A_k, B_{j,k}, N(k)$ be the same as in {\bf Case 2} and $M(u)=\max(k\in \N: d_k\leq \delta_2(u))$. Clearly
$$\bigcup_{k=0}^{M(u)-1}A_k\subset[0,\delta_2(u)]\subset\bigcup_{k=0}^{M(u)}A_k,\ \ \ \bigcup_{k=0}^{N(k)-1}B_{j,k}\subset A_k\subset\bigcup_{k=0}^{N(k)}B_{j,k},$$
and by Bonferroni's inequality
\BQN\label{3lower}
\pi(u)&\geq&\sum_{k=0}^{M(u)-1}\sum_{j=0}^{N(k)-1}\pk{\sup_{t\in B_{j,k}}X(t)>u}-
\underset{(j,k)\prec(j',k')}{\sum_{(j,k),(j',k')\in\mathcal{L}'}}\pk{\sup_{t\in B_{j,k}}X(t)>u,\sup_{t\in B_{j',k'}}X(t)>u}\nonumber\\
&=:& J'_1(u)-J'_2(u),
\EQN
where $\mathcal{L}'=\{(j,k):0\leq k\leq M(u)-1, 0\leq j\leq N(k)-1\}$.\\
By \eqref{sig}, \nelem{lein}, \nelem{leYZ} and similar argumentation as \eqref{J1} with $\mathcal{G}_u^{\pm\vn}$ replaced by $\mathcal{F}_u^{\pm\vn}$ and the fact that
$(\ln u)^{1/\beta}d_{M(u)+1}\rw \IF$,  $u\rw\IF$, we get
\BQN\label{3lower1}
\lim_{S\rw\IF}\lim_{u\rw\IF}\frac{J'_1(u)(\ln u)^{1/\beta}}{u^{2/\alpha}\Psi(u)}
\geq a^{1/\alpha}\mathcal{H}_{\alpha}\int_{0}^{\IF}e^{\frac{-2b x^\beta}{\alpha^2}}dx.
\EQN
By\cite{atlocally}[Lemma 4.5]
\BQN\label{3lower2}
J'_2(u)&=&\underset{(j,k)\prec(j',k')}{\sum_{(j,k),(j',k')\in\mathcal{L}'}}\pk{\sup_{t\in B_{j,k}}X(t)>u,\sup_{t\in B_{j',k'}}X(t)>u}\nonumber\\
&\leq&\underset{(j,k)\prec(j',k')}{\sum_{(j,k),(j',k')\in\mathcal{L}'}}\pk{\sup_{t\in B_{j,k}}\overline{X}(t)>u,\sup_{t\in B_{j',k'}}\overline{X}(t)>u}\nonumber\\
&=&o\LT(u^{2/\alpha}(\ln u)^{-1/\beta}\Psi(u)\RT),\ u\rw\IF.
\EQN
Hence inserting \eqref{3lower1} and \eqref{3lower2} into \eqref{3lower}, we have
\BQNY
\lim_{u\rw\IF}\frac{\pi(u)(\ln u)^{1/\beta}}{u^{2/\alpha}\Psi(u)}\geq a^{1/\alpha}\mathcal{H}_{\alpha}\int_{0}^{\IF}e^{\frac{-2b x^\beta}{\alpha^2}}dx,
\EQNY
which combined with \eqref{3upper} gives that
\BQNY
\pi(u)\sim a^{1/\alpha}\mathcal{H}_{\alpha}u^{2/\alpha}(\ln u)^{-1/\beta}\Psi(u)\int_{0}^{\IF}e^{\frac{-2b x^\beta}{\alpha^2}}dx,
\ u\rw\IF.
\EQNY
Consequently, according to \nelem{le1} and
\BQNY
\pi(u)\leq\pk{\sup_{t\in[0,T]}X(t)>u}\leq \pi(u)+\pk{\sup_{t\in[\delta(u),T]}X(t)>u},
\EQNY
\eqref{nee} is proved and all claims follow.
\QED
\section{Appendix}
In this section we present the proofs of the lemmas used in the proof of \netheo{Main1}. \\
\prooflem{le1}
Below $\mathbb{Q}_k,\ k=0,1,2\ldots $, are some positive constants.\\
{\bf Step 1}: First we prove \eqref{le1eq1}.
By the continuity of $\sigma(t)$ in [0,T], for any small enough constant $0<\theta<1$
\BQNY
\sup_{t\in[\theta,T]}\sigma(t)=:\rho(\theta)<\sigma(t_0)=\sigma(0)=1.
\EQNY
Then by Borell inequality in \cite{AdlerTaylor}
\BQNY
\pk{\sup_{t\in[\theta,T]}X(t)>u}\leq \exp\LT(-\frac{\LT(u-\mathbb{Q}_0\RT)^2}{2\rho^2(\theta)}\RT)=o\LT(\Psi\LT(u\RT)\RT),
\EQNY
as $u\rw\IF$, where $\mathbb{Q}_0=\E{\sup_{t\in[0,T]}X(t)}<\IF$.\\
By assumption (iv), for any small $\vn\in(0,1)$, when $\theta$ small enough
\BQNY
1+(1-\vn)ce^{-|t|^{-\gamma}}\leq\frac{1}{\sigma(t)}\leq1+(1+\vn)ce^{-|t|^{-\gamma}},
\EQNY
holds for $t\in [0,\theta]$.
Then
\BQNY
\frac{1}{\sigma(t)}\geq 1+(1-\vn)ce^{-|t|^{-\gamma}}\geq 1+(1-\vn)cu^{-2}(\ln u)^q
\EQNY
uniformly holds for $t\in [\delta_1(u),\theta]$.\\
Moreover by assumption (i) and (iii), when $\theta$ small enough
\BQNY
\E{(X(t)-X(s))^2}&=&\E{X^2(t)}+\E{X^2(s)}-2\E{X(t)X(s)}\\
&\leq&2-2(1-2a(t)|t-s|^{\alpha(t)})\\
&\leq&\mathbb{Q}_1|t-s|^{\varsigma}
\EQNY
holds uniformly for $s,t\in [0,\theta]$, where $\mathbb{Q}_1=\sup_{t\in[0,\theta]}4a(t)$ and  $\varsigma=\inf_{t\in[0,\theta]}\alpha(t)>0$.\\
Then by Piterbarg inequality
\BQNY
\pk{\sup_{t\in[\delta_1(u),\theta]}X(t)>u}\leq \mathbb{Q}_2\theta u^{2/\varsigma}
\Psi(u[1+(1-\vn)cu^{-2}(\ln u)^q])=o\LT(\Psi\LT(u\RT)\RT), \ \ u\rw\IF.
\EQNY
Further, since
\BQNY
\pk{\sup_{t\in[0,T]}X(t)>u}\leq\pk{\sup_{t\in[0,\delta_1(u)]}X(t)>u}+\pk{\sup_{t\in[\delta_1(u),\theta]}X(t)>u}
+\pk{\sup_{t\in[\theta,T]}X(t)>u},
\EQNY
and
\BQNY
\pk{\sup_{t\in[0,T]}X(t)>u}\geq\pk{\sup_{t\in[0,\delta_1(u)]}X(t)>u}\geq \pk{X(0)>u}=\Psi\LT(u\RT),
\EQNY
we get
\BQNY
\pk{\sup_{t\in[0,T]}X(t)>u}\sim\pk{\sup_{t\in[0,\delta_1(u)]}X(t)>u}, \ \ u\rw\IF.
\EQNY
{\bf Step 2}: Next we prove \eqref{le1eq2}.
When $\gamma\leq \beta$, since $\delta_1(u)=o(\delta_2(u))$, as $u\rw\IF$ and
by {\bf Step 1}
\BQNY
\pk{\sup_{t\in[\delta_1(u),T]}X(t)>u}=o\LT(\Psi\LT(u\RT)\RT), \ \ u\rw\IF.
\EQNY
Then for $u$ large enough, \eqref{le1eq2} is obvious.\\
When $\gamma>\beta$, for $u$ large enough,  we have $\delta_2(u)<\delta_1(u)$ and
\BQNY
\pk{\sup_{t\in[\delta_2(u),T]}X(t)>u}\leq \pk{\sup_{t\in[\delta_2(u),\delta_1(u)]}X(t)>u} +\pk{\sup_{t\in[\delta_1(u),T]}X(t)>u}.
\EQNY
By {\bf Step 1}, we know for all $u$ large
\BQNY
\pk{\sup_{t\in[\delta_1(u),T]}X(t)>u}\leq \Psi(u),
\EQNY
and then we just need to deal with $\pk{\sup_{t\in[\delta_2(u),\delta_1(u)]}X(t)>u}$.

Since $\delta_1(u)\rw 0,\ u\rw\IF$,  then by assumption (v)
\BQNY
\alpha(t)>\alpha+\frac{3}{4}b(\delta_2(u))^\beta
\EQNY
holds for all $t\in[\delta_2(u),\delta_1(u)]$ when $u$ large enough.\\
Let $\eta_u=u^{-2/\LT(\alpha+\frac{3}{4}b(\delta_2(u))^\beta\RT)}$. For sufficiently large $u$ and $s,t\in[\delta_2(u),\delta_1(u)]$, there exists a constant $\mathbb{Q}_3>0$ such that
\BQNY
1-r(s,t)\leq1-e^{-\mathbb{Q}_3|s-t|^{\alpha+\frac{3}{4}b(\delta_2(u))^\beta}}.
\EQNY
Let $Y_u(t), t\geq 0$ be a family of centered stationary Gaussian processes with correlation functions
\BQNY
r_Y(s,t)=e^{\mathbb{Q}_3|s-t|^{\alpha+\frac{3}{4}b(\delta_2(u))^\beta}}.
\EQNY
Then from Slepian's inequality we get for any constant $S>0$
\BQNY
\pk{\sup_{t\in[\delta_2(u),\delta_1(u)]}X(t)>u}&\leq& \pk{\sup_{t\in[\delta_2(u),\delta_1(u)]}\frac{X(t)}{\sigma(t)}>u}\\
&\leq& \pk{\sup_{t\in[\delta_2(u),\delta_1(u)]}Y_u(t)>u}\\
&\leq& \pk{\sup_{t\in[0,S]}Y_u(t)>u}\\
&\leq& \sum_{i=0}^{\lfloor S\eta_u^{-1}\rfloor+1} \pk{\sup_{t\in[i\eta_u,(i+1)\eta_u]}Y_u(t)>u}\\
&\leq& (\lfloor S\eta_u^{-1}\rfloor+1) \pk{\sup_{t\in[0,\eta_u]}Y_u(t)>u},
\EQNY
for sufficiently large $u$. Notice that for each $s,t\in[0,1]$
\BQNY
1-r_Y(\eta_u t,\eta_u s)=\mathbb{Q}_3u^{-2}|s-t|^{\alpha+\frac{3}{4}b(\delta_2(u))^\beta}(1+o(1))
=\mathbb{Q}_3u^{-2}|s-t|^{\alpha}(1+o(1)),\ u\rw\IF.
\EQNY
Hence, from \cite{Pit96}[Lemma D.1]
\BQNY
 \pk{\sup_{t\in[0,\eta_u]}Y_u(t)>u}\sim\mathcal{H}_\alpha[1]\Psi(u),
\EQNY
as $u\rw\IF$. Combining this with the fact that
\BQNY
\eta_u^{-1}&=&u^{2/\LT(\alpha+\frac{3}{4}\delta_2(u)\RT)}=u^{2/\alpha}u^{2/\LT(\alpha+\frac{3}{4}\delta_2(u)\RT)-2/\alpha}
=u^{2/\alpha}u^{-\frac{3}{2}(\delta_2(u))^\beta/\LT(\alpha\LT(\alpha+\frac{3}{4}(\delta_2(u))^\beta\RT)\RT)}\\
&=&u^{2/\alpha}u^{-\frac{3}{2}\frac{\alpha^2 (\ln (\ln u))}{\beta(\ln u)}/\LT(\alpha\LT(\alpha+\frac{3}{4}(\delta_2(u))^\beta\RT)\RT)}
\leq u^{2/\alpha}u^{-\frac{4}{3}\frac{\ln (\ln u)}{\beta(\ln u)}}= u^{2/\alpha}(\ln u)^{-4/(3\beta)},
\EQNY
we get for some constant $\mathbb{Q}_4$ and all $u$ large enough
\BQNY
 \pk{\sup_{t\in[\delta_2(u),\delta_1(u)]}X(t)>u}\leq \mathbb{Q}_4 S u^{2/\alpha} (\ln u)^{-4/3\beta}\Psi\LT(u\RT).
\EQNY
Then the result follows.
\QED

\BEL\label{lein}
 Under the notation in the proof of \netheo{Main1}, for $(j,k)\in \mathcal{U}=\{(j,k):0\leq k\leq M^*(u), 0\leq j\leq N(k)\}$ and $\lim_{u\rw\IF}\frac{f(u)}{u}=1$, there exists $u_0$ such that for each $u\geq u_0$\\
\underline{1)} $\pk{\sup_{t\in B_{j,k}}\overline{X}(t)>f(u)}\geq\pk{\sup_{t\in [0,S]}Y_{\nu,u}(t)>f(u)}$;\\
\underline{2)} $\pk{\sup_{t\in B_{j,k}}\overline{X}(t)>f(u)}\leq\pk{\sup_{t\in [0,S]}Z_{\nu,u}(t)>f(u)},$\\
where
\BQNY
M^*(u)=\LT\{
\begin{array}{ll}
0,&\ \hbox{if}\ \gamma<\beta,\\
M_\en(u),&\ \hbox{if}\ \gamma=\beta,\\
M(u),&\ \hbox{if}\ \gamma>\beta.\\
\end{array}
\RT.
\EQNY
\EEL
\prooflem{lein}
Since the proofs of scenarios $\gamma<\beta,\ \gamma=\beta,$ and $\gamma>\beta$ are similar, we only present the proof of $\gamma=\beta$.
Set $X_{j,k,u}(t)=\overline{X}\LT(d_k+\frac{jS+t}{u^{2/\alpha(d_k)}}\RT)$, then  $ \sup_{t\in B_{j,k}}\overline{X}(t) \overset{d}{=} \sup_{t\in [0,S]}X_{j,k,u}(t).$
It is enough to analyze the supremum of $X_{j,k,u}(t)$.\\
\COM{The rest of the proof is similar  as the proof of  \cite{atlocally} [Lemma 4.1] where mainly use the Slepian's inequality, if we notice to replace $t_u$ in [Lemma 4.1] with $\delta(u)$.\\}
\underline{1)} For sufficiently large $u$ and $s,t\in[0,T]$
\BQN\label{lein1}
1-Cov\LT(X_{j,k,u}(s),X_{j,k,u}(t)\RT)&=&1-Cov\LT(\overline{X}\LT(d_k+\frac{jS+s}{u^{2/\alpha(d_k)}}\RT),
\overline{X}\LT(d_k+\frac{jS+t}{}\RT)\RT)\nonumber\\
&\geq&(1-\nu/2)^{1/3}a\LT|u^{-2/\alpha(d_k)}(s-t)\RT|^{\alpha\LT(d_k+u^{-2/\alpha(d_k)}(jS+t)\RT)}\nonumber\\
&=&(1-\nu/2)^{1/3}a u^{-2\alpha\LT(d_k+u^{-2/\alpha(d_k)}(jS+t)\RT)/\alpha(d_k)}\LT|(s-t)\RT|^{\alpha\LT(d_k+u^{-2/\alpha(d_k)}(jS+t)\RT)}\nonumber\\
&=&(1-\nu/2)^{1/3}a \times I_1\times I_2.
\EQN
We deal with $I_1$ and $I_2$ separately.
For sufficiently large $u$, uniformly with respect to $k$,
\BQN\label{lein2}
I_1&=&u^{-2\alpha\LT(d_k+u^{-2/\alpha(d_k)}(jS+t)\RT)/\alpha(d_k)}\nonumber\\
&=&u^{-2}u^{2\LT(\alpha(d_{k})-\alpha\LT(d_k+u^{-2/\alpha(d_k)}(jS+t)\RT)\RT)/\alpha(d_k)}\nonumber\\
&=&u^{-2}e^{2(\ln u)\LT(\alpha(d_{k})-\alpha\LT(d_k+u^{-2/\alpha(d_k)}(jS+t)\RT)\RT)/\alpha(d_k)}\nonumber\\
&\geq&u^{-2}(1-\nu/2)^{1/3},
\EQN
where the last inequality follows from the fact that
\BQNY
(\ln u)\LT|\alpha(d_{k})-\alpha\LT(d_k+u^{-2/\alpha(d_k)}(jS+t)\RT)\RT|
&\leq&(\ln u)\LT(\LT|b(d_{k})^\beta-b\LT(d_k+u^{-2/\alpha(d_k)}(jS+t)\RT)^\beta\RT|+2\delta_1^{\beta+\delta}(u)\RT)\\
&\leq&(\ln u)\LT(\frac{b}{(\ln u)(\ln \ln u)^{1/\beta}}+2\delta_1^{\beta+\delta}(u)\RT)\\
&\leq&\frac{b}{(\ln \ln u)^{1/\beta}}+2(\ln u)\LT(\frac{1}{2\ln u-q\ln\ln u}\RT)^{\frac{\beta+\delta}{\gamma}}\ \rw 0, \ u\rw\IF.
\EQNY
For $I_2$, we need to prove that
\BQN\label{I2}
I_2\geq (1-\nu/2)^{1/3}|s-t|^{\alpha+2b\delta_1^\beta(u)}.
\EQN
Assumption (v) implies that
\BQN\label{I22}
\alpha\LT(d_k+u^{-2/\alpha(d_k)}(jS+t)\RT)<\alpha+2b\delta_1^\beta(u)
\EQN
for each $(j,k)\in\mathcal{U}$.  Thus if $|s-t|<1$, then \eqref{I2} holds immediately.
If $1\leq|s-t|\leq S$, then by \eqref{I22}
\BQNY
I_2&=&\LT|(s-t)\RT|^{\alpha\LT(d_k+u^{-2/\alpha(d_k)}(jS+t)\RT)}\\
&\geq&T^{\alpha\LT(d_k+u^{-2/\alpha(d_k)}(jS+t)\RT)-\alpha-2b\delta_1^\beta(u)}|s-t|^{\alpha+2b\delta_1^\beta(u)}\\
&\geq&T^{-2b\delta_1^\beta(u)}|s-t|^{\alpha+2b\delta_1^\beta(u)}\\
&\geq&(1-\nu/2)^{1/3}|s-t|^{\alpha+2b\delta_1^\beta(u)}
\EQNY
for sufficiently large $u$.
The above combined with \eqref{lein1}, \eqref{lein2} and \eqref{I2} gives that for sufficiently large $u$, uniformly with respect to
$(j,k)\in \mathcal{U}$,
\BQNY
1-Cov\LT(X_{j,k,u}(s),X_{j,k,u}(t)\RT)\geq (1-\nu/2) au^{-2}|s-t|^{\alpha+2b\delta_1^\beta(u)}\geq 1-Cov\LT(Y_{\nu,u}(s),Y_{\nu,u}(t)\RT).
\EQNY
Thus by  Slepian's inequality \underline{1)} is proved.\\
\underline{2)}
For all $u$ large
\BQNY
1-Cov\LT(X_{j,k,u}(s),X_{j,k,u}(t)\RT)&=&1-Cov\LT(\overline{X}\LT(d_k+\frac{jS+s}{u^{2/\alpha(d_k)}}\RT),
\overline{X}\LT(d_k+\frac{jS+t}{}\RT)\RT)\nonumber\\
&\leq&(1+\nu)^{1/3}a\LT|u^{-2/\alpha(d_k)}(s-t)\RT|^{\alpha\LT(d_k+u^{-2/\alpha(d_k)}(jS+t)\RT)}.
\EQNY
Following the argument analogous to that for the proof of \underline{1)}, we obtain that for sufficiently large $u$, uniformly with respect to $k$, and $s,t\in[0,S]$
\BQNY
1-Cov\LT(X_{j,k,u}(s),X_{j,k,u}(t)\RT)\leq 1-Cov\LT(Z_{\nu,u}(s),Z_{\nu,u}(t)\RT).
\EQNY
Again the application of Slepian's inequality  completes the proof. \QED

\BEL\label{leYZ}
 For $S>1$, $\nu\in(0,1)$, and $\lim_{u\rw\IF}\frac{f(u)}{u}=1$, as $u\rw\IF$, we have\\
\underline{1)} $\pk{\sup_{t\in [0,S]}Y_{\nu,u}(t)>f(u)}=\mathcal{H}_{\alpha}\LT[0,S((1-\nu)a)^{1/\alpha}\RT]\Psi\LT(f(u)\RT)(1+o(1))$;\\
\underline{2)} $\pk{\sup_{t\in [0,S]}Z_{\nu,u}(t)>f(u)}=\mathcal{H}_{\alpha}\LT[0,S((1+\nu)a)^{1/\alpha}\RT]\Psi\LT(f(u)\RT)(1+o(1)).$
\EEL
\prooflem{leYZ} We present the proof of \underline{1)} and omit the proof of  \underline{2)} since it follows with similar arguments.
Following the definition of $Y_{\nu,u}(t)$, for each $s,t\in[0,S]$
\BQNY
&&\lim_{u\rw\IF}f^2(u)\LT[1-Cov\LT(Y_{\nu,u}\LT(t(a(1-\nu))^{-1/\alpha}\RT),Y_{\nu,u}\LT(s(a(1-\nu))^{-1/\alpha}\RT)\RT)\RT]\\
&&\quad\quad=\lim_{u\rw\IF}\LT(a(1-\nu)\RT)^{1-\LT(\alpha+2b\delta^\beta(u)\RT)/\alpha}|s-t|^{\alpha+2b\delta^\beta(u)}=|s-t|^\alpha.
\EQNY
Moreover, for all $s,t\in[0,S]$, sufficiently large $u$ and some constant $C>0$
\BQNY
&&f^2(u)\LT[1-Cov\LT(Y_{\nu,u}\LT(t(a(1-\nu))^{-1/\alpha}\RT),Y_{\nu,u}\LT(s(a(1-\nu))^{-1/\alpha}\RT)\RT)\RT]\\
&&\quad\quad \leq\LT(a(1-\nu)\RT)^{1-\LT(\alpha+2b\delta^\beta(u)\RT)/\alpha}|s-t|^{\alpha+2b\delta^\beta(u)}\leq CT^{2\alpha}|s-t|^\alpha,
\EQNY
where the last inequality follows from the fact that
\BQNY
|s-t|^{\alpha+2b\delta^\beta(u)}\leq |s-t|^\alpha, \ \hbox{if}\ |s-t|<1,
\EQNY
and
\BQNY
|s-t|^{\alpha+2b\delta^\beta(u)}\leq T^{2\alpha}\leq  T^{2\alpha}|s-t|^\alpha, \ \hbox{if}\ 1\leq|s-t|\leq T.
\EQNY
Hence, by \cite{HP2004}[Lemma 7], we conclude that
\BQNY
\pk{\sup_{t\in [0,S]}Y_{\nu,u}(t)>f(u)}
&=&\pk{\sup_{t\in [0,((1-\nu)a)^{1/\alpha}S]}Y_{\nu,u}((a(1-\nu))^{-1/\alpha}t)>f(u)}\\
&=&\mathcal{H}_{\alpha}\LT[0,((1-\nu)a)^{1/\alpha}S\RT]\Psi\LT(f(u)\RT)(1+o(1)),
\EQNY
as $u\rw\IF$. This completes the proof. \QED

{\bf Acknowledgement}: Thanks to  Swiss National Science Foundation grant no.  200021-166274.

\bibliographystyle{plain}
\bibliography{atmodelBBC}
\end{document}